\documentstyle[12pt]{article}
\catcode`\@=11
\@addtoreset{equation}{section}

\catcode`\@=12
\newtheorem{Theorem}{Theorem}[section]
\newtheorem{Definition}{Definition}[section]
\newtheorem{Proposition}{Proposition}[section]
\newtheorem{Lemma}{Lemma}[section]
\newtheorem{Corollary}{Corollary}[section]

\title{Pre-Lagrangian Submanifolds
In Contact Manifolds\thanks{Project 19871044 Supported by NSF}}

\author{Renyi Ma\\
Department of Applied Mathmatics \\
Tsinghua University \\
Beijing, 100084\\
People's Republic of China}

\date { }

\begin{document}
\textwidth=165mm
\textheight=230mm
\parindent=8mm
\frenchspacing
\maketitle

\begin{abstract}
In this article, we prove the non-existence of 
exact pre-Lagrangian submanifolds in 
contact manifolds by using the Gromov's 
nonlinear Fredholm alternative. 
\end{abstract}

\noindent{\bf 1991 MR Subject Classification} 58F, 58G, 53C

\noindent{\bf Chinese Library Clasification} O175.2, O186.11

\section{Introduction and Results}

Let $U$ be a $(2n-1)$-dimensional manifolds. A contact 
structure $\xi $ on $U$ is a completely nonintegrable codimension 
$1$ tangent distribution. It means that $\xi $ can be defined, 
at least locally, by a $1-$form 
$\lambda $ with $\lambda \wedge (d\lambda )^{n-1}\ne 0$. 
Note that if $n$ is odd then the contact distribution 
$\xi $ is automatically orientable. For an even $n$ the existence 
of a contact structure implies the orientability of the ambient manifold 
$U$. In both cases, the coorientability of 
$\xi $ implies that $\xi $ and $U$ are both orientable. 
We will asume from now on that $\xi $ is coorientable and fix its orientation. 
Then $\xi $ can be globally defined by a $1-$form $\lambda $, 
which is determined up to a multiplication by a positive 
function. Let $SU=U\times ]0,\infty [$.  We will 
still denote by $\lambda $ the pull back of $\lambda $ by the projection 
of $SU=U\times R_{+}$ on the first factor 
and denote by $t$ the projection on the second. 
Then the form $\omega =d(a\lambda )$ defines a symplectic structure 
on $SU$(indeed, 
$(d(a\lambda ))^n=a^{n-1}da\wedge \lambda \wedge d\lambda 
^{n-1}\ne 0)$. The map $(x,a)\to (x,a/f(x))$ induces an isomorphism 
of forms $d(a\lambda )$ and $d(a\mu)$ for $\mu =f\lambda $. Therefore, 
the symplectic manifold $(SU, \omega )$ depends, up 
to a symplectomorphism, only on the contact manifolds $(U, \xi )$ and 
not on the choice of the $1-$form $\lambda $. For an $n-$dimensional 
manifolds $M$ let us denote by $P_+T^*(M)$ the oriented projective cotangent 
bundle of $M$ with the contact structure 
$\xi $ defined by the 
form $pdq$. The manifold $P_+T^*M$ can also be considered 
a space of cooriented $(n-1)$-dimensional contact elements of 
$M$. With this interpretation the plane $\xi _x$ of $\xi $ at a 
point $x=(p,q)$, $q\in M$, $p\in T^*_q(M)$, consists 
of infinitesimal deformations of $\xi _x$, which 
leaves fixed the point of contact $q\in M$. Then the symplectization 
$Sympl(P_+T^*(M),\xi )$ is isomorphic 
to $T^*(M)\setminus M$ with the standard symplectic structure 
$\omega =d(pdq )$, for more example, see[1-6].

According to [2,4], the following notion
was suggested by D. Bennequin. 

An $n-$dimensional submanifolds $L$ of the $(2n-1)-$dimensional 
contact manifold $(U, \xi )$ is called $pre-Lagrangian$ if 
it satisfies the following two conditions:

a.  $L$ is transverse to $\xi $;

b.  The distribution $\xi \cap T(L)$ is integrable and 
can be defined by a closed $1-$form. 

\vskip 3pt 

For any pre-Lagrangian submanifold $L\subset U$ there 
exists a Lagrangian submanifold $\hat {L}\subset S_{\xi }U$ 
such that $\pi (\hat {L})=L$. The cohomology 
class $\lambda \in H^1(L;R)$, such that 
$\pi ^*\lambda =[\alpha _\xi |\hat {L}]$, is 
defined uniquely up 
to multiplication by a non-zero constant. Conversly, 
if $L\subset U$ is the 
(embedded) image of a Lagrangian submanifold 
$\hat {L}\subset S_{\xi }U$ under 
the projection $S_{\xi }U\to U$ then 
$L$ is pre-Lagrangian(see[2,4]). Thus 
with any pre-Lagrangian submanifold 
$L\subset U$ one can canonically associate a projective 
class of the form $\lambda $. The main result 
of this paper is following:

\begin{Theorem}
There does not exist any pre-Lagrangian 
submanifold $L\subset U$ with the canonical projective class 
equal to zero, especially any simply connected 
manifold can not be embedded in $(U,\xi )$ as 
a pre-Lagrangian submanifold. 
\end{Theorem}

\begin{Theorem}
Let $(U,\lambda )$ be 
a closed contact manifold and $\varphi : L 
\to (U, \lambda )$ a closed Pre-Lagrangian embedding, 
then $[\varphi ^*(\lambda )]\ne 0$ in
$H^*(L,R)$, especially
$H^1(L)\ne 0$.
\end{Theorem}

{\bf Sketch of proofs}: We will work in the framework 
proposed by Gromov in [5]. In Section 2, we study the 
linear Cauchy-Riemann operator and sketch some basic 
properties. In section 3, we study 
the space ${\cal D}(V,W)$
of contractible disks in manifold $V$ with boundary 
in Lagrangian submanifold $W$ and construct a Fredholm
section of tangent bundle of  ${\cal  D }(V,W)$.  
In Section 4, we
use the Gromov's trick in [5] to 
estimate the energy of the solutions 
of the nonlinear Cauchy-Riemann equations.
In the final section, we
use Gromov's nonlinear Fredholm trick to 
complete our proof as in [5].

\section{Linear Fredholm Theory}

For $100<k<\infty $ consider the Hilbert space
$V_k$ consisting of all maps $u\in H^{k,2}(D, C^n)$,
such that $u(z)\in R^n\subset C^n$ for almost all $z\in
\partial D$. $L_{k-1}$ denotes the usual Hilbert $L_{k-1}-$space
$H_{k-1}(D, C^n)$. We define an operator
$\bar \partial :V_p\mapsto L_p$ by
\begin{equation}
\bar \partial u=u_s+iu_t
\end{equation}
where the coordinates on $D$ are
$(s,t)=s+it$, $D=\{ z||z|\leq 1\}  $.
The following result is well known(see[5]).
\begin{Proposition}
$\bar \partial :V_p\mapsto L_p $ is a surjective real
linear Fredholm operator of index $n$. The kernel
consists of the constant real valued maps.
\end{Proposition}
Let $(C^n, \sigma =-Im(\cdot ,\cdot ))$ be the standard
symplectic space. We consider a real $n-$dimensional plane
$R^n\subset C^n$. It is called Lagrangian if the
skew-scalar product of any two vectors of $R^n$ equals zero.
For example, the plane $p=0$ and $q=0$ are Lagrangian
subspaces. The manifold of all (nonoriented) Lagrangian subspaces of
$R^{2n}$ is called the Lagrangian-Grassmanian $\Lambda (n)$.
One can prove that the fundamental group of
$\Lambda (n)$ is free cyclic, i.e.
$\pi _1(\Lambda (n))=Z$. Next assume
$(\Gamma (z))_{z\in \partial D}$ is a smooth map
associating to a point $z\in \partial D$ a Lagrangian
subspace $\Gamma (z)$ of $C^n$, i.e.
$(\Gamma (z))_{z\in \partial D}$ defines a smooth curve
$\alpha $ in the Lagrangian-Grassmanian manifold $\Lambda (n)$.
Since $\pi _1(\Lambda (n))=Z$, one have
$[\alpha ]=ke$, we call integer $k$ the Maslov index
of curve $\alpha $ and denote it by $m(\Gamma )$, see([1]).

Now let $z:S^1\mapsto R^n\subset C^n$ be a smooth curve. 
Then
it defines a constant loop $\alpha $ in Lagrangian-Grassmanian
manifold $\Lambda (n)$. This loop defines
the Maslov index $m(\alpha )$ of the map
$z$ which is easily seen to be zero.

  Now Let $(V,\omega )$ be a symplectic manifold and 
$W\subset V$ a closed Lagrangian submanifold. Let 
$u:D^2\to V$ be a smooth map homotopic to constant map 
with boundary 
$\partial D\subset W$.
Then $u^*TV$ is a symplectic vector bundle and 
$(u|_{\partial D})^*TW$ be a Lagrangian subbundle in 
$u^*TV$. Since $u$ is contractible, we can take 
a trivialization of $u^*TV$ as 
$$\Phi (u^*TV)=D\times C ^n$$ 
and 
$$\Phi (u|_{\partial D})^*TW)\subset S^1\times C ^n$$
Let 
$$\pi _2: D\times C ^n\to C^n$$
then 
$$\bar u: z\in S^1\to \pi _2\Phi (u|_{\partial D})^*TW(z)\in \Lambda (n).$$
Write $\bar u=u|_{\partial D}$.
\begin{Lemma}
Let $u: (D^2,\partial D^2) \rightarrow (V,W)$ be a $C^k-$map $(k\geq 1)$ 
as above. Then,
$$m(u|_{\partial D})=0$$
\end{Lemma}
Proof.  Since $u$ is contractible in $V$ relative to $W$, we have 
a homotopy $\Phi _s$ of trivializations
such that 
$$\Phi _s(u^*TV)=D\times C ^n$$ 
and 
$$\Phi _s(u|_{\partial D})^*TW)\subset S^1\times C ^n$$
Moreover 
$$\Phi _0(u|_{\partial D})^*TW=S^1\times R^n$$
So, the homotopy induces 
a homotopy $\bar h$ in Lagrangian-Grassmanian
manifold. Note that $m(\bar h(0, \cdot ))=0$.
By the homotopy invariance of Maslov index,
we know that $m(u|_{\partial D})=0$.

\vskip 5pt 

   Consider the partial differential equation
\begin{eqnarray}
\bar \partial u+A(z)u=0  \ on \ D  \\
u(z)\in \Gamma (z) R^n\ for \ z\in \partial D \\
\Gamma (z)\in GL(2n,R)\cap Sp(2n)\\
m(\Gamma )=0 \ \ \ \ \ \ \ \ 
\end{eqnarray}

For $100<k<\infty $ consider the Banach space $\bar V_k $
consisting of all maps $u\in H^{k,2}(D, C^n)$ such 
that  $u(z)\in \Gamma (z)$ for almost all $z\in
\partial D$. Let $L_{k-1}$ the usual $L_{k-1}-$space $H_{k-1}(D,C^n)$ and 

$$L_{k-1}(S^1)=\{ u\in H^{k-1}(S^1)|u(z)\in \Gamma (z) R^n\ for \ z\in 
\partial D\}$$ 
We define an operator $P$:
$\bar V_{k}\rightarrow L_{k-1}\times L_{k-1}(S^1)$ by
\begin{equation}
P(u)=(\bar \partial u+Au,u|_{\partial D})
\end{equation}
where $D$ as in (2.1).
\begin{Proposition}
$\bar \partial : \bar V_p \rightarrow L_p$
is a real linear Fredholm operator of index n.
\end{Proposition}

\section{Nonlinear Fredholm Theory}

\noindent{\bf 3.1. Adapted metrics in symplectic manifold 
$(M,\omega )$.} A Riemannian 
metric $g$ on $M$ is called adapted (to the 
symplectic form $\omega $) if 
$g+\sqrt{-1}\omega $ is a Hermitian metric 
with respect to some almost complex 
structure $J:T(M)\to T(M)$ preserving $g$ and $\omega $. 
This is equivalent to the existence of a 
$g-$orthonormal coframe $x_i,y_i$,$i=1,...,n=dimM/2$, at each 
point in $M$ such that $\omega $ equals 
$\sum _1^nx_i\wedge y_i$ at this point. Yet another 
equivalent definition reads 
$$||dH||_g=||grad_{\omega }H||_g$$
for all smooth functions $H$ on $M$, where, recall, $grad_{\omega }H$ 
is the (Hamiltonian) 
vector field which is $\omega -$dual to $dH$. 

Let us show that a complete adapted metric always exists.
\begin{Lemma}
(Eliashberg-Gromov[3]). Every symplectic manifold
$M=(M,\omega )$ admits a complete adapted metric $g$.
\end{Lemma}
Proof(due to [3]). The required metric will be constructed starting 
with arbitrary adapted metric $g_0$ and applying a certain 
symplectic automorphism $A$ of 
$T^*(M)$ to it. This $A$ is constructed with an exhaustion of 
$M$ by compact domains 
with smooth boundaries $S_i$ expands $g_0$ transversally 
to all $S_i$. Namely, we take small 
$\varepsilon _i-$neighbourhoods $N_i\subset M$ of $S_i$, normally 
(with respect to $g_0$) decomposed as $N_i=S_i\times [-\varepsilon _i,
\varepsilon _i]$. We denote by $\Sigma _i\subset T(N_i)$ and 
$\nu _i\in T(N_i)$ the subbundles tangent and 
normal to the slices $S_i\times t$, $t\in [-\varepsilon _i,
\varepsilon _i]$, respectively and take some 
symplectic automorphisms $A_i:T(N_i)\to T(N_i)$ preserving the decomposition 
$T(N_i)=\Sigma _i\oplus \nu _i$ and acting on $\nu _i$ by 
$A_i(\nu )=2\nu $. Then $A$ 
is taken equal to $Id$ outside all 
$N_i$ and $A|T(N_i)=^{def}A^{\varphi _i}$ where 
$\varphi _i(s,t)=\varphi _i(t)$ is a suitable 
sequence of positive functions on $[-\varepsilon _i,\varepsilon _i]$ 
such that $\varphi _i$ vanish at ends $\pm \varepsilon _i$ and are 
large and fast growing with $i$ on the subsegments 
$[-\varepsilon _i/2,\varepsilon _i/2]$. Clearly 
$g=Ag_0$ is complete(as well as adapted) for suitable $\varphi _i$. 

\vskip 3pt

\noindent{\bf 3.2. Construction of Lagrangian submanifolds.} 
Let $(V',\omega ')$($\omega =d\alpha '$) be an exact 
symplectic manifold and $W'\subset V'$ a closed 
submanifolds, we call $W'$ an exact Lagrangian submanifold if 
$\alpha '|W'$ an exact form, i.e., 
$\alpha '|W'=df$. Consider an isotopy of Lagrange submanifolds 
in $V'$ given by a $C^{\infty }-$map $F':W'\times [0,1]\to V'$ and 
let $\tilde {\omega }'$ be the pull-back of the form $\omega '$ 
to $W'\times [0,1]$. The form $\tilde {\omega }'$ clearly is exact 
since $\omega '=d\alpha '$, 
$\tilde {\omega }'=d\tilde {l}'$, where the $1-$form $\tilde {l}'$ 
is closed on $W'\times t$ for $t\in [0,1]$. Recall that $F'$ is called an 
$exact $ $isotopy$ if the class 
$[\tilde {l}'|W'\times t]\in H^1(W'=W'\times t;R)$ is constant in 
$t\in [0,1]$, for more detail see[$5,2.3B'$].

\vskip 3pt

Let $U$ a contact manifold and $L\subset U$ an exact 
pre-Lagrangian submanifold, as proved in [4], that 
one can choose a contact form $\lambda $ on $U$ 
such that 
$(V',\omega ')=(U\times R_+,d(a\lambda ))$ and 
$d(a\lambda )|L\times \{1\}=0$ and $\lambda |L$ is 
exact. So $W'(=\{1\}\times L)\subset V'$ an 
exact Lagrangian submanifold in $V'$
and the manifold $(SU,d(a\lambda ))$ has a canonical 
diffeotopy $v'\to sv'$ for $s\in [0,\infty )$. The induced isotopy on 
$L$ clearly is Lagrange; it is exact if and only if the form $l'|W'$ 
is exact. The isotopied manifolds $W'_s=s(W')$ are disjointed from 
$W'$ for any $s$.
We choose a positive number $s_0$ small enough which 
will be determined in section 5 and 
define 
$$F':W'\times [0,1]\to SU$$
as 
$$F'((w,1),t)=(w,1+ts_0)$$
Then one can easily check that $F'$ is 
an exact Lagrangian isotopy of $W'$ in $SU$.  

\vskip 2pt 

Let 
$(V,\omega )=(V'\times C,\omega '\oplus \omega _0)$. 
As in [5], we use figure eight trick to 
construct a Lagrangian submanifold in $V$ through the 
Lagrange isotopy $F'$ in $V'$. 
Fix a positive $\delta <1$ and take a $C^{\infty }$-map $\rho :S^1\to 
[0,1]$, where the circle $S^1$ ia parametrized by $\Theta \in [-1,1]$, 
such that the $\delta -$neighborhood $I_0$ of $0\in S^1$ goes to 
$0\in [0,1]$ and $\delta -$neighbourhood $I_1$ of $\pm 1\in S^1$ 
goes $1\in [0,1]$. Let 
\begin{eqnarray}
\tilde {l}&=&-\psi (w',\rho (\Theta ))\rho '(\Theta )d\Theta \cr
&=&-\Phi d\Theta  
\end{eqnarray}
be the pull-back of the form 
$\tilde {l}'=-\psi (w',t)dt $ to $W'\times S^1$ under the map 
$(w',\Theta )\to (w',\rho (\Theta ))$ and 
assume without loss of generality $\Phi $
vanishes on $W'\times (I_0\cup I_1)$. 

  Next, consider a map $\alpha $ of the annulus $S^1\times [\Phi _-,\Phi _+]$ 
into $R^2$, where $\Phi _-$ and $\Phi _+$ are the lower and the upper 
bound of the fuction $\Phi $ correspondingly, such that 
   
   $(i)$ The pull-back under $\alpha $ of the form 
$dx\wedge dy$ on $R^2$ equals $-d\Phi \wedge d\Theta $. 
  
   $(ii)$ The map $\alpha $ is bijective on $I\times [\Phi _-,\Phi _+]$ 
where $I\subset S^1$ is some closed subset, 
such that $I\cup I_0\cup I_1=S^1$; furthermore, the origin 
$0\in R^2$ is a unique double point of the map $\alpha $ on 
$S^1\times 0$, that is 
$$0=\alpha (0,0)=\alpha (\pm 1,0),$$  
and 
$\alpha $ is injective on $S^1=S^1\times 0$ minus $\{ 0,\pm 1\}$. 

   $(iii)$ The curve $S^1_0=\alpha (S^1\times 0)\subset R^2$ ``bounds'' 
zero area in $R^2$, that is $\int _{S^1_0}xdy=0$, for the $1-$form 
$xdy$ on $R^2$. 
\begin{Proposition}
Let $V'$, $W'$ and $F'$ as above. Then there exists  
an exact Lagrangian embedding $F:W'\times S^1\to V'\times R^2$ 
given by $F(w',\Theta )=(F'(w',\rho (\Theta )),\alpha (\Theta ,\Phi ))$.
\end{Proposition}
Proof. Similar to [5,2.3$B_3'$].

\vskip 3pt 

\noindent{\bf 3.3. Formulation of Hilbert manifolds}.  
Now let $(U,\lambda )$ be a contact
manifold with contact form
$\lambda$. Let $SU=(U\times ]0,\infty [,d(a\lambda )$
be its symplectization. 
By Lemma 3.1, one has 

\begin{Proposition}
There exists an adapted complete metric on 
the symplectization $SU=(U\times ]0,\infty [,d(a\lambda ))$ 
of contact manifolds $(U,\lambda )$.
\end{Proposition}
In the following we denote by 
$(V,\omega )=(SU\times R^2,d(a\lambda )\oplus dx\wedge dy))$ 
with the adapted metric $g\oplus g_0 $ and 
$W\subset V$($W=F(W'\times S^1)$) the Lagrangian submanifold constructed 
in section 3.2.

\vskip 3pt 

   Let $k\geq 100$ and 
$${\cal D}^k(V,W,p)=\{ u \in H^k(D,V)|
u(\partial D)\subset W,\ u\ 
homotopic \ to \ u_0=p, \ u(1)=p\}$$
\begin{Lemma}
Let $W$ be a closed Lagrangian submanifold in 
$V$. Then, 
${\cal D}^k(V,W,p)$
is a Hilbert manifold with the tangent bundle
\begin{equation}
T{\cal D}^k(V,W,p)=\bigcup _{u\in {\cal {D}}^k(V,W,p)}
\Lambda ^{k-1}(u^*TV,u|_{\partial D}^*TW,p)
\end{equation}
here 
\begin{eqnarray}
\Lambda ^{k-1}(u^*TV,u|_{\partial D}^*TW,p)=\ \ \ \ \ \ \\ 
\{ H^{k-1}-sections \ of \ (u^*(TV),(u|_{\partial D})^*TL)\ 
which \ vanishes  \ at \ 1\} 
\end{eqnarray}
\end{Lemma}
Proof: See[5].

   Now we construct a nonlinear Fredholm operator
from ${\cal D}^k(V,W,p)$ to $T{\cal D}^k(V,W, p)$ follows in [5].
Let $\bar \partial :{\cal D}^k(V,W,p)\rightarrow T{\cal D}^k(V,W,p)$
be the Cauchy-Riemmann Section induced by the Cauchy-Riemann
operator, locally,
\begin{equation}
\bar \partial u={{\partial u}\over {\partial s}}
+J{{\partial u}\over {\partial t}}
\end{equation}
for $u\in {\cal D}^k(V,W,p)$.

Since the space ${\cal D}^k(V,W,p)$ is Hilbert manifold, the
tangent space $T{\cal D}^k(V,W,p)$ is trivial, i.e.
there exists a bundle isomorphism
$$\Phi : T{\cal {D}}^k(V,W, p)
\rightarrow {\cal {D}}^k(V,W,p)
\times E$$
where $E$ is a Hilbert Space.
Then the Cauchy-Riemann section $\bar \partial $
on $T{\cal D}^k(V,W,p) $ induces a nonlinear map
$$\Phi \circ \bar \partial : {\cal D}^k(V,W, p)\mapsto E$$
In the following, we still denote
$\Phi \circ \bar \partial $ by $\bar \partial $
for convenience.
Now we define
\begin{eqnarray}
F: {\cal D}^k(V,W,p)\rightarrow E\\
F(u)=\Phi (\bar \partial u)
\end{eqnarray}

\begin{Theorem}
The nonlinear operator $F$ defined in (3.6-3.7)
is a nonlinear Fredholm operator of Index zero.
\end{Theorem}
Proof. According to the definition of the nonlinear
Fredholm operator, we need to prove that
$u\in {\cal D}^k(V,W,p)$, the linearization
$DF(u)$ of $F$ at $u$ is a linear Fredholm
operator.
Note that
\begin{equation}
DF(u)=D{\bar \partial _{[u]}}
\end{equation}
where
\begin{equation}
(D\bar \partial _{[u]})v=\frac{\partial v}{\partial s}
+J\frac{\partial v}{\partial t}+A(u)v
\end{equation}
with 
$$v|_{\partial D}\in (u|_{\partial D})^*TW$$
here $A(u)$ is $2n\times 2n$
matrix induced by the torsion of
almost complex structure, see [5] for the computation.

   Observe that the linearization $DF(u)$ of 
$F$ at $u$ is equivalent to the following Lagrangian 
boundary value problem
\begin{eqnarray}
{{\partial v}\over {\partial s}}+J{{\partial v}\over {\partial t}}
+A(u)v=f, \ v\in \Lambda ^k(u^*TV)\\ 
v(t)\in T_{u(t)}W, \ \ t\in {\partial D}
\end{eqnarray}
One 
can check that (3.10-11) 
defines a linear Fredholm operator. In fact, 
by proposition 2.2 and Lemma 2.1, since the operator $A(u)$ is a compact, 
we know that the operator $F$ is a nonlinear Fredholm operator 
of the index zero.

\begin{Definition}
A nonlinear Fredholm $F:X\rightarrow Y$ operator is
proper if any $y\in Y$, $F^{-1}(y)$ is finite or for
any compact set $K\subset Y$, $F^{-1}(K)$ is compact
in $X$.
\end{Definition}
\begin{Definition}
$deg(F,y)=\sharp \{ F^{-1}(y)\} mod2$ is called the Fredholm
degree of a nonlinear proper Fredholm operator(see[5,11]).
\end{Definition}
\begin{Theorem}
Assum that the nonlinear Fredholm operator
$F: {\cal D}^k(V,W,p)\rightarrow E$
constructed in (3.6-7) is proper. Then,
$$deg(F,0)=1$$
\end{Theorem}
Proof: We assume that $u:D\mapsto V$ be a $J-$holomorphic disk
with boundary $u(\partial D)\subset W$. Since almost complex
structure $\widetilde {J}$ tamed by  the symplectic form $\omega $,
by stokes formula, we conclude $u: D^2\rightarrow
w$ is a constant map. Because $u(1)=p$, We know that
$F^{-1}(0)={p}$.
Next we show that the linearizatioon $DF(p)$ of $F$ at $p$ is
an isomorphism from $T^p{\cal D}(V,W,p)$ to $E$.
This is equivalent to solve the equations
\begin{eqnarray}
{\frac {\partial v}{\partial s}}+J{\frac {\partial v}{\partial t}}=f\\
v|_{\partial D}\subset T_pW
\end{eqnarray}
By Lemma 3.1, we know that $DF(p)$ is an isomorphism.
Therefore $deg(F,0)=1$.
\begin{Corollary}
$deg(F,w)=1$ for any $w\in E$.
\end{Corollary}
Proof. Using the connectedness of $E$ and the
homotopy invariance of $deg$.

\section{Non-properness of Fredholm Operator }

We shall prove in this section that the operator
$F:{\cal D}\rightarrow E$ constructed
in the above section is non proper along the line in [5].

\vskip 3pt 

\noindent{\bf 4.1. Anti-holomorphic section.} 
  Let $C=R^2$ and $(V',\omega ')$,  
$(V,\omega )=(V'\times C, \omega '\oplus \omega _0)$, 
and 
$W$ as in section 3 and $J=J'\oplus i$, $g=g'\oplus g_0$, 
$g_0$ the standard metric on $C$. 

   Now let $c\in C$(here $C$ the complex plane) 
be a non-zero vector. We consider the 
equations
\begin{eqnarray}
v=(v',f):D\to V'\times C \nonumber \\
\bar \partial _{J'}v'=0,\bar \partial f=c\ \ \ \nonumber \\
v|_{\partial D}:\partial D\to W\ \ \ 
\end{eqnarray}
here $v$ homotopic to constant map 
$\{ p\}$ relative to $W$. 
Note that $W\subset V\times B_R(0)$(here $R$ 
depends on the $s_0$ in section 3.2). 
\begin{Lemma}
Let $v$ be the solutions of (4.1), then one has 
the following estimates
\begin{eqnarray}
E({v})=
\{ 
\int _D(g'({{\partial {v'}}\over {\partial x}},
{J'}{{\partial {v'}}\over {\partial x}})
+g'({{\partial {v'}}\over {\partial y}},
{J'}{{\partial {v'}}\over {\partial y}}) \nonumber \\
+g_0({{\partial {f}}\over {\partial x}},
{i}{{\partial {f}}\over {\partial x}})
+g_0({{\partial {f}}\over {\partial y}},
{i}{{\partial {f}}\over {\partial y}}))d\sigma \}
\leq 4\pi R^2. 
\end{eqnarray}
\end{Lemma}
Proof: Since $v(z)=(v'(z),f(z))$ satisfy (4.1)
and $v(z)=(v'(z),f(z))\in V'\times C$ 
is homotopic to constant map $v_0:D\to \{ p\}\subset W$ 
in $(V,W)$, by the Stokes formula
\begin{equation}
\int _{D}v^*(\omega '\oplus \omega _0)=0
\end{equation}
Note that the metric $g$ is adapted to the symplectic form 
$\omega $ and $J$, i.e., 
\begin{equation}
g=\omega  (\cdot ,J\cdot )
\end{equation}
By the simple algebraic computation, we have 
\begin{equation}
\int _{D}{v}^*\omega  ={{1}\over {4}}
\int _{D^2}(|\partial v|^2 
-|\bar {\partial }v|^2)=0
\end{equation}
and 
\begin{equation}
|\nabla v|={{1}\over {2}}(
|\partial v|^2 +|\bar \partial v|^2 
\end{equation}
Then 
\begin{eqnarray}
E(v)&=&\int _{D} |\nabla v| \nonumber \\ 
      &=&\int _{D}\{ {{1}\over {2}}(
|\partial v|^2+|\bar \partial v|^2)\} d\sigma \nonumber \\ 
&=&\pi |c|_{g_0}^2
\end{eqnarray}
By the equations (4.1), 
one get 
\begin{equation}
\bar \partial f=c \ on \ D
\end{equation}
We have 
\begin{equation}
f(z)={{1}\over {2}}c\bar z+h(z)
\end{equation}
here $h(z)$ is a holomorphic function on $D$. Note that  
$f(z)$ is smooth up to the boundary $\partial D$, then, by 
Cauchy integral formula
\begin{eqnarray}
\int _{\partial D}f(z)dz&=&{{1}\over {2}}c\int _{\partial D}
\bar {z}dz+\int _{\partial D}h(z)dz \cr
&=&\pi ic
\end{eqnarray}
So, we have 
\begin{equation}
|c|={{1}\over {\pi}}|\int _{\partial D^2}f(z)dz|
\end{equation}
Therefore, 
\begin{eqnarray}
E(v)&\leq &\pi |c|^2
\leq {{1}\over {\pi }}|\int _{\partial D}f(z)dz|^2      \cr
&\leq &{{1}\over {\pi }}|\int _{\partial D}|f(z)||dz||^2   \cr
&\leq &4\pi |diam(pr_2(W))^2 \cr
&\leq &4\pi R^2.
\end{eqnarray}
This finishes the proof of Lemma.

\begin{Proposition}
For $|c|\geq 3R$, then the 
equations (4.1)
has no solutions. 
\end{Proposition}
Proof. By (4.11), we have 
\begin{eqnarray} 
|c|&\leq &{{1}\over {\pi }}\int _{\partial D}|f(z)||dz|\cr
&\leq &{{1}\over {\pi }}\int _{\partial D}
diam(pr_2(W))||dz| \cr
&\leq &2R
\end{eqnarray}
It follows that $c=3R$ can not be obtained by 
any solutions.

\vskip 3pt 

\noindent{\bf 4.2. Modification of section $c$.} 
Note that the section $c$ is not a section of the 
Hilbert bundle in section 3 since $c$ is not 
tangent to the Lagrangian submanifold $W$, we must modify it as follows:

\vskip 3pt 

  Let $c$ as in section 4.1, we define 
\begin{eqnarray}
c_{\chi ,\delta }(z,v)=\left\{ \begin{array}{ll}
c \ \ \ &\mbox{if\  $|z|\leq 1-2\delta $,}\cr
0 \ \ \ &\mbox{otherwise}
\end{array}
\right. 
\end{eqnarray}
Then by using the cut off function $\varphi _h(z)$ and 
its convolution with section 
$c_{\chi ,\delta }$, we obtain a smooth section 
$c_\delta$ satisfying

\begin{eqnarray}
c_{\delta }(z,v)=\left\{ \begin{array}{ll}
c \ \ \ &\mbox{if\  $|z|\leq 1-3\delta $,}\cr
0 \ \ \ &\mbox{if\  $|z|\geq 1-\delta $.}
\end{array}
\right. 
\end{eqnarray}
for $h$ small enough by well-known convolution theory.

   Now let $c\in C$ be a non-zero vector and 
$c_\delta $ the induced anti-holomorphic section. We consider the 
equations
\begin{eqnarray}
v=(v',f):D\to V'\times C \nonumber \\
\bar \partial _{J'}v'=0,\bar \partial f=c_\delta \ \ \ \nonumber \\
v|_{\partial D}:\partial D\to W\ \ \ 
\end{eqnarray}
Note that $W\subset V\times B_{R}(0)$ for $2\pi R^2$. 
Then by repeating the same argument as section 4.1., we obtain 
\begin{Lemma}
Let $v$ be the solutions of (4.16) and $\delta $ 
small enough, then one has 
the following estimates
\begin{eqnarray}
E({v})\leq 4\pi R^2. 
\end{eqnarray}
\end{Lemma}
and

\begin{Proposition}
For $|c|\geq 3R$, then the 
equations (4.16)
has no solutions. 
\end{Proposition}

\begin{Theorem}
The Fredholm operator $F: {\cal  {D}}^k(V,W,p) \rightarrow E$ is not proper.
\end{Theorem}
Proof. If $F$ is proper, taking 
a path $\gamma (\mu )$ connecting 
$0$ and $c$, then $F^{-1}(\gamma (\cdot ))$ 
is a compact set in 
${\cal D}^k(V,W,p)$ for $k\geq 100$, then 
the gradients of map $v$ have a uniforms 
bounds, i.e., 
\begin{equation}
|\nabla v|\leq c_1 \ \ for \ v\in F^{-1}(\gamma (\cdot )) 
\end{equation}
Note that $v(1)=p$, the above bounds imply 
\begin{equation}
v|_{\partial D}\subset W'(0)\times (-K,K)
\end{equation}
for $K$ large enough which only depends on $c_1$. 
Since $W'(0)\times (-k,+k)$ is regular embedding in 
$V$, we know that 
$v|_{\partial D}$ is compact in the submanifold
$W$. Then Theorem 3.1 and 3.2, i.e., that the index of $F$
is zero and $deg(F)=1$ implies $F$ can take 
the value $c$ for $c\geq 3R$, This 
contradicts Proposition 4.1.
So, $F$ is not proper.

\section{Nonlinear Fredholm Alternative}

In this section, we use the Sacks-Uhlenbeck-Gromov's trick and 
Gromov's nonlinear Fredholm alternative to prove
the existence of $J$-holomorphic disk with 
boundary in $W$ if $W\subset SU\times C$ is 
Lagrangian submanifold. 

\

{\bf Proof of Theorem 1.1}. If Theorem 1.1 does not hold, i.e., 
there exists a exact Pre-Lagrangian submanifold $L$ 
in a contact manifold $U$, 
we use the canonical isotopy in the symplectization 
to construct an very small Lagrangian 
isotopy of $L$ then by the Gromov's figure eight construction 
in section 3.2 we obtain the exact Lagrangian submanifold 
$W$ in $SU\times C$. By choosing $s_0$ in section 3.2 small 
enough such that $4\pi R^2$ small enough we conclude that 
the solutions of (4.16) is bounded by using the 
monotone inequality of minimal surface since the boundary of solutions 
of (4.16) remain in the compact manifold $W$. Then 
for large vector $c\in C$ in equations (4.16) we know
that the nonlinear 
Cauchy-Riemann equations  
has no solution, this implies that the operator $F$ constructed 
in section 3.3 is not proper or the solutions of equations (4.16) is non
-compact. 
The non-properness of the 
operator implies 

\vskip 2pt 

{\bf a.} The existence of $J-$holomorphic plane $v:C\to V$ with 
bounded energy $E(v)\leq E_0$. Since $v$ has 
a bounded image then by Gromov's removal singularity theorem 
we get a non constant map $w:S^2\to V$ which 
contradict the exactness of $V$.

\vskip 2pt 

{\bf b.} The existence of $J-$holomorphic half plane 
$v:H\to V$ with boundary $\partial H$ 
in $W$. Since $v$ has a bounded image, then by the Gromov's 
removal boundary singularity we get a 
$J-$holomorphic disks $w:D\to V$ with boundary 
in $W$, this contradicts that $W$ is an exact Lagrangian 
submanifold. 

\vskip 3pt 

This implies 
Theorem 1.1 holds.

\end{document}